# NEURAL NETWORK AND SEGMENTED LABOUR MARKET[1]

Patrice GAUBERT, METIS-MATISSE, Université de Paris I
Marie COTTRELL, SAMOS-MATISSE, Université de Paris I

## 1. Introduction

In France, for administrative reasons, unemployed workers may actually be involved in occasional work while remaining identified as unemployed (and receiving the corresponding benefit). This is due to the fact that the unemployed are deemed to be seeking full-time jobs and non-fixed term contracts of employment. This situation may be analysed as evidence of a special type of secondary segment of the labour market in a context of massive unemployment. The authors consider the effects of this situation both on the duration of unemployment and its recurrence may be usefully investigated.

This paper is a contribution to the analysis of segmentation in the French labour market. It is based on a very different approach from the existing literature which is focused on career differentiation (Favereau, 1991; Theodossiou, 1995), returns to education (Balsan, 1994), or, more generally, efficiency wage (Oi,1990 ; Albrecht, 1992) or differentiation and discrimination (Dickens, 1985 ; Boston, 1990).
Basically labour market segmentation is described as comprising a primary and a secondary segment, differentiated mainly by job types : permanent in the primary segment, precarious in the secondary one.
In contrast to the standard presentation, the high level of unemployment observed in France and the growth of precarious employment (part-time work, fixed-term contracts of employment, etc.), raises questions about the relationship between the secondary segment and unemployment. The rapid growth of these new job categories, often low skilled and low paid jobs, may be interpreted in two ways .Either they illustrate a change in the nature of the worker's transition from unemployment to employment (easing the search for a new job by a newly unemployed individual thanks to the closening of his links with the labour market), or, on the contrary, are they evidence of a rapid extension of the secondary segment of the labour market coupled with a symmetrical shrinkage of the primary segment, one of the characteristic features of recent economic growth trends.
High unemployment in France raises questions as to what extent this model can be applied to the structure of unemployment. Is it possible to differentiate among categories of unemployed workers according to the duration and recurrence of unemployment ? How can these categories be defined, using information about educational and skill level, work experience and other aspects of previous jobs and the reasons for their termination ?

[1] *This study was conducted with financial assistance from the Government Employment Service (ANPE). The authors accept full responsability for the ideas expressed herewithin.*
*The authors wish to thank Patrick Rousset, participants of seminar SAMOS at Université de Paris I, and two referees for valuable suggestions and Leslie Thompson (University of Paris I) for his invaluable assistance in drafting the English language version of this article.*

These questions are analysed using a Kohonen network, to classify the unemployed.
The data used are drawn from the register of unemployed kept by the Government Employment Service (ANPE). Two regions chosen were Nord and Rhone-Alpes. Some forty thousands individuals were randomly selected from the register in each of these two regions between July 1993 and August 1996.
Our results show that unemployed workers who accept occasional work experience a significantly longer period of unemployment. For these workers the longer the time spent doing casual work during this period, the shorter the overall duration of the latter. The corresponding figure seems to be higher than that corresponding to the total monthly hours occupied in occasional work. The likelihood of recurrence of unemployment appears greater for such individuals than for other unemployed workers, while the length of their unemployment is close to that measured for people who do no occasional work.

**2. Data**
The complete register of the unemployed held by the ANPE contains information on all registered unemployed workers in the twenty-six administrative regions of France. The sample used for this study is made up of individuals who were looking for work between 1 July 1993 and 31 August 1996 and who were still registered as seeking a job or having found one by the latter date.
As this data was originally collected for administrative purposes, it had to be adjust to be used here.
One of the main problems we were faced with is that an individual may be registered more than once. This is what is meant by recurrence of unemployment. Each entry records the particular circumstances and category of job termination, how unemployment ended, duration of unemployment and the time spent in occasional work. This information is essential for determining the factors correlated to the duration of unemployment.
Our investigation of individual workers revealed that some of them had been registered more than once and so this data was collated. This was possible when dealing with quantitative data (length of each period of unemployment, amount of months in occasional work, number of job offers received through the Government Employment Service). We also recorded the different periods of unemployment. But the qualitative characteristics of each search for work could not be taken into consideration. Only those concerning the latest search were available anyway. Individual characteristics are assumed to remain the same over this relatively short period of time.
The proportion of individuals registered only once in two regions studied is recorded in the table below (the figures are similar for the other administrative regions) :

Table 1

| Region | Total nb of records | Total nb of individuals | Individuals with one registration | % (3/2) |
|---|---|---|---|---|

|  | (1) | (2) | (3) |  |
|---|---|---|---|---|
| Rhône-Alpes | 1 310 887 | 883 385 | 585 153 | .66 |
| Nord | 1 185 835 | 760 348 | 475 520 | .62 |

### 3. Variables created concerning occasional work

The data concerning casual work reveals some problems of measurement. The figures for total monthly hours are far from accurate, because only those concerning the first and last months of casual work are known and the exact dates of beginning and end are not available either. So a change in the total monthly hours worked from one month to the next one may mean a job change or that a job started for instance the twentieth of the previous month or ended the eighth of the following month. Other information relating to other records is not known precisely either. For example, data is available only concerning the range within which the real measure is located (less than or greater than seventy-eight hours per month). The wage data for casual work is also sometimes missing.

For each period of unemployment involving occasional work we made a number of approximations. The missing information was replaced by computing the mean monthly hours for individuals belonging to the same administrative class. For each individual we then compared the average figure for the full duration of the occasional work with the maximum monthly hours declared during this period. As the difference is less than 10 % we decided to use the maximum observed in each period, because this figure actually exists and this avoids the problem of incomplete months. Similarly we used the maximum wage declared, which we assumed to be correctly associated with the monthly hours as defined.

Different variables were used in classifying and interpreting the results. First we used variables describing the individual (age, marital status, number of children, nationality) and professional characteristics (educational and skill level, seniority in the latest job). Second, we chose to examine different aspects of the worker's unemployment status, such as the reason for his job loss (lay off, voluntary quit, end of fixed term contract, first job search, return to labour market after a temporary exit), how his/her unemployment status ended (job, training program, exit from labour market, cancellation of registration for administrative reasons), duration and cumulated duration of unemployment measured in days and months (for individuals with more than one unemployment period), number of job offers received through the Employment Service. We also investigated the nature of the occasional work (total monthly hours, number of months in this work and cumulated duration, number of different periods of occasional work within a given period of unemployment, monthly wage).

In order to avoid level effects we use relative values to determine the impact of occasional work on the duration of unemployment: the total months of occasional work divided by the total months or total years of unemployment, total periods of occasional work divided by the total number of unemployment periods and by total years of unemployment, the number of job offers divided by total months of unemployment.

The figures obtained (expressed in means) are reported in the tables below.

Table 2

| Region | Obs | Duration | Cumulated duration | Nb of periods | Nb of job | Job/Unemp. |
|---|---|---|---|---|---|---|
| NORD | 38068 | 476.52 | 636.44 | 1.55 | 0.50 | 0.11 |
| RHONE | 44165 | 410 | 547 | 1.48 | 0.52 | 0.09 |

| Region | Job/year | Offers/Unemp. | Wages | Hours | Age | Seniority |
|---|---|---|---|---|---|---|
| NORD | 0.54 | 0.12 | 767.08 | 37.17 | 31.32 | 4.32 |
| RHONE | 0.50 | 0.11 | 1345.68 | 33.45 | 32.45 | 4.65 |

On the whole the Nord and Rhone regions are quite similar in terms of the variables presented, except, as far as the length of unemployment is concerned, which is slightly greater for the Nord region (twenty-one months versus eighteen).

The distribution of this sample is reported in table 3, 4 and 5 by using frequencies computed with different variables to be developed later in our analysis.

Table 3

| | Duration | | | Number of periods | | | Age bracket | | |
|---|---|---|---|---|---|---|---|---|---|
| | < 6 m. | 6-12 m. | > 12 m. | 1 p. | 2 p. | > 2 p. | < 35 | 35 à 55 | > 55 |
| NORD | 37.3 | 20.8 | 41.9 | 62.8 | 25.0 | 12.2 | 66.5 | 30.7 | 2.8 |
| RHONE | 39.7 | 21.7 | 38.6 | 66.1 | 23.1 | 10.8 | 64.3 | 31.7 | 4.0 |

| | Skill | | | Level of education | | | | | |
|---|---|---|---|---|---|---|---|---|---|
| | Sup. | Clerks | Indus | S4 | S2 | Sec dip. | Second. | Prim. | None |
| NORD | 10.9 | 49.5 | 39.6 | 4.9 | 8.0 | 14.5 | 39.9 | 7.0 | 25.7 |
| RHONE | 17.5 | 51.0 | 31.5 | 8.8 | 11.3 | 14.5 | 39.7 | 6.2 | 19.5 |

Sup. : executive, technicians
Indus : skilled and unskilled workers
S4 : Master level and more
S2 : Undergraduate level
Sec. Dip. : Secondary level with high school diploma

Table 4

| | Reasons for job loss | | | | | Type of unemployment exit | | | | |
|---|---|---|---|---|---|---|---|---|---|---|
| | l. o. | end | 1st job | return | na | job | p. job | quit | training | cancel | f. job |
| NORD | 15.1 | 54.5 | 16.9 | 10.3 | 3.2 | 25.8 | 13.2 | 10.5 | 6.8 | 13.4 | 30.3 |
| RHONE | 23.9 | 46.7 | 11.2 | 12.7 | 5.5 | 27.8 | 18.8 | 9.6 | 6.4 | 14.1 | 23.3 |

l.o. : lay-off
end : end of fixed-term contract
p.job : registration cancelled by administration, most of the time because job found
f.job : individuals occupying a full-time job while seeking another one (not unemployed)

Table 5

|  | Monthly hours in occasional work | | | | | Nb of casual job per period | | | Proportion of casual work during period of unemployment | | | |
|---|---|---|---|---|---|---|---|---|---|---|---|---|
|  | 0 | 0-39 | 39-78 | 78-117 | > 117 | 0 | 0-1 | > 1 | 0 | 0-0.1 | 0.1-0.3 | > 0.3 |
| NORD | 67.6 | 2.6 | 7.1 | 3.8 | 18.9 | 59.4 | 29.4 | 11.2 | 59.4 | 11.9 | 14.2 | 14.5 |
| RHONE | 67.4 | 2.6 | 9.6 | 7.1 | 13.3 | 62.4 | 25.3 | 12.3 | 62.4 | 11.5 | 14.9 | 11.3 |

Some regional differences should however be pointed out. These relate to educational levels (higher levels were found in the Rhone region) and to the relative shares of lay-offs, ends of contract and first jobs.

Data concerning the unemployment allowance received by unemployed people was not available, though this is an important variable as far as our study is concerned. Information on the kind of job found after unemployment was not available either. We were thus not able to carry out a detailed analysis of the duration of unemployment and of the overall situation of the unemployed embracing consecutive periods of work and unemployment. Hence our study will be limited to identifying the main categories of unemployed workers.

### 3. Classification and class characteristics

The Kohonen classification shown in table 6 (Kohonen, 1984, 1985 ; Cottrell, 1997, *see in the Appendix I for some details about the method*) is obtained using quantitative variables defined for each individual as age, job seniority, length of the latest period of unemployment, cumulated duration of unemployment since initial registration, number of job offers received per month of unemployment, number of periods of unemployment, the duration of occasional jobs per month of unemployment, number of periods spent doing occasional work per year of unemployment and per number of periods of unemployment, and finally, total monthly hours worked and total wages received for occasional work.

The variables defined earlier to describe occasional work are at the core of our analysis. The number of periods spent doing occasional work compared separately to the number of periods of unemployment *and* to the duration of this unemployment is used to differentiate between people presumably engaged in a work close to what they had doing before losing their job (occasional work occupying a large part of the period of unemployment and only one period of occasional work in each period of unemployment) and those people engaged in very occasional work for short periods of time and who change occupations within the same period of unemployment (presumably with no relation between these jobs and the one lost due to unemployment).

Two maps are presented in the Appendix II, showing the code vector of each class (the line joining the different values is not intended to plot an evolution : the values are those obtained at the end of the algorithm for normalised variables ranked in alphabetical order according to name). These maps also include the ten super-classes resulting from a clustering classification applied to the Kohonen classification.

Table 6 - Kohonen classification (super-classes)

| SUPER-CLASS | Size | Age | Seniority | Duration last period | Cumulated duration | Number of periods |
|---|---|---|---|---|---|---|
| **NORD** | | | | | | |
| 1 | 9472 | 24.23 | 1.30 | 183.49 | 241.24 | 1.26 |
| 2 | 4790 | 29.75 | 2.19 | 485.12 | 930.82 | 1.77 |
| 3 | 3047 | 47.01 | 22.05 | 440.69 | 539.29 | 1.28 |
| 4 | 2745 | 45.00 | 3.59 | 320.31 | 376.18 | 1.17 |
| 5 | 2954 | 33.25 | 4.55 | 1072.32 | 1159.83 | 1.16 |
| 6 | 1563 | 28.33 | 4.20 | 168.05 | 236.03 | 1.41 |
| 7 | 3139 | 28.20 | 2.62 | 164.90 | 639.94 | 3.57 |
| 8 | 3678 | 37.31 | 4.56 | 1617.14 | 1673.05 | 1.03 |
| 9 | 4610 | 27.57 | 2.76 | 354.32 | 525.66 | 1.77 |
| 10 | 2070 | 28.03 | 3.39 | 158.16 | 206.79 | 1.32 |
| Total | 38068 | 31.32 | 4.32 | 476.52 | 636.44 | 1.55 |
| **RHONE** | | | | | | |
| 1 | 10239 | 25.28 | 1.50 | 261 | 264 | 1.02 |
| 2 | 2808 | 25.75 | 1.62 | 212 | 325 | 1.51 |
| 3 | 8447 | 28.89 | 2.38 | 191 | 623 | 2.62 |
| 4 | 3949 | 42.91 | 4.82 | 220 | 304 | 1.20 |
| 5 | 5194 | 48.57 | 19.34 | 446 | 545 | 1.27 |
| 6 | 3524 | 36.48 | 2.30 | 1380 | 1415 | 1.05 |
| 7 | 2691 | 31.20 | 3.24 | 462 | 725 | 1.80 |
| 8 | 2200 | 28.76 | 4.81 | 140 | 182 | 1.29 |
| 9 | 2870 | 31.94 | 2.74 | 1069 | 1118 | 1.09 |
| 10 | 2243 | 30.58 | 4.72 | 257 | 307 | 1.30 |
| Total | 44165 | 32.45 | 4.65 | 410 | 547 | 1.48 |

The above table and those figuring in the appendix II show the main characteristics of the super-classes in terms of the means computed for the variables used in the classification. Our findings are significant in both regions studied.

Two groups were found to be made up of older people. One group was characterised by lengthy seniority (*Nord* 3, *Rhone* 5) while seniority in the other was found to be relatively short (*Nord* 4, *Rhone* 4). This may be due to recent changes in the career of workers in the first group (unemployment caused by obsolete skills or because they were working in greatly depressed industrial sectors within the region) while in the second group workers experience long periods of unemployment soon after their entry into the labour market.

Two groups comprising people of intermediate age (*Nord* 5 and 8, *Rhone* 6 and 9) were characterised by a single lengthy period of unemployment. It was observed that some people have already experienced a relatively long period of unemployment at the date of registration.

Other groups point to differences between the two regions, such as systematic recurrence of unemployment in all other groups in the Nord region, while this phenomenon is present only in some groups in the Rhone region. The reason for this lies in the high percentage of « ends of contract » in the Nord region among the causes of unemployment.

The main characteristics for the Nord region can be summarised as follows. One group is characterised by a high rate of recurrence (7) and there are two groups in which recurrence is a significant feature (2 and 9). One other group is both particularly young and numerous (1). Two other groups are characterised by a rather short period of unemployment (6 and 10). As for the Rhone region two groups show a high degree of recurrence (3 and 7), and two others show a particularly short length of unemployment (8 and 10). Another group is young and has experienced a single short period of unemployment (1), quite close to group 2 in age but different from it due to regular recurrence.

The data furnished by the distribution of the population within these groups, using both initial qualitative variables and those constructed using ranges of quantitative variables, has proven invaluable in enabling a more accurate understanding of the different categories of unemployed workers. It should be pointed that in both regions there is one group (*Nord* 8 and *Rhone* 10) that has to be distinguished from the rest of the population studied. This group contains a high proportion of people corresponding to an administrative category of workers seeking work though already occupying a full-time job. This category is not deemed acceptable in so far as the French administration mixes these people with workers doing occasional work.

Another result, which lends support to our earlier interpretation obtained with quantitative variables, is that the two groups composing the first category (*Nord* 3 and *Rhone* 5) are principally made up of laid-off workers possessing high skills but low levels of education.

Groups 2, 7 and 9 for the Nord region and group 3 for the Rhone region are characterised by a high percentage of « ends of contract », while groups 7 and 3 were also mentioned above to have employment profiles characterised by short duration. This would seem to provide evidence of a mechanism linking precarious work and unemployment. It can be seen that for some groups finding a job is an important reason for leaving unemployment (*Nord* 6, 9 and 10, *Rhone* 1, 2 and 8) and this is in turn linked to a shorter duration of unemployment (*Nord* 6, 7 and 10, *Rhone* 3 and 8). Doing occasional work also plays an important role for some of these groups (*Nord* 5, 9 and 10, *Rhone* 2, 7 and 9).

We will now turn to a closer examination of the different classes obtained.

### 4. Class interpretation

To obtain a more precise description of these classes and of what characterises the individuals belonging to them we carried out a multiple correspondence analysis using some of the qualitative variables already presented or which we created using some of the initial quantitative ones. We added the classification variable representing the super-class in which each individual is placed. We selected the following ten variables, in addition to the one furnishing the super-class number (the values plotted on the graphs figuring in the appendix II are in brackets)

- **number of children** : no child (*0enf*), one child or more (*enf*)
- **educational level** : secondary school not completed (*form0*), secondary school completed (*nbacc*), post secondary school level (*+bacc*)
- **skills** : various classes of workers defined according to their skills (*mano*, *ousp*, *op12*, *o3hq*), unskilled clerks (*emnq*), skilled clerks (*empq*), technicians (*tecd*), foremen and executives (*mait*, *cadr*)
- **reasons for unemployment** : end of contract (*fin*), lay-off (*lic*), first job (*prem*), return to labour market (*repr*)
- **duration in months of latest period of unemployment** : less than 6 months (*inf6*), between 6 and 12 months (*inf12*), more than 12 months (*sup12*)
- **type of exit from unemployment**: job (*empl*), mostly job (*pemp*), training (*stag*), departure from labour market (*retr*), cancellation for administrative reasons (*susp*)
- **recurrence of unemployment** : single period of unemployment (*rec1*), two periods (*rec2*), three periods or more (*rec3*)
- **average monthly hours of occasional work during the latest period of unemployment** : no occasional work (*0ar*), between 0 and 39 hours (*arm39*), from 40 to 78 hours (*ar78*), from 79 to 117 hours (*ar117*), more than 117 hours (*arp117*)
- **percentage of cumulated duration of unemployment doing occasional work** : no work (*0par*), below 10 % (*01par*), between 10 and 30 % (*13par*), above 30 % (*p3par*)
- **average number of periods of occasional work per year of unemployment** : no work (*0peran*), less than one period (*m1peran*), one or two periods (*12peran*), more than two periods (*p2peran*).

Some of these variables (duration in months of latest period of unemployment, recurrence of unemployment over the whole period, percentage of cumulated duration in occasional work, monthly hours spent in these jobs during the latest period of unemployment) are obtained by discretization of the quantitative variables used to produce the classification. Others were already qualitative ones.

These eleven variables together with their fifty-five subcategories produce a complete cross tabulation table and a Burt table (extension of a contingency table when there are more than two variables). We applied a multiple correspondence analysis to this table producing successive planes with decreasing informational content representing simultaneously the fifty-five subcategories. Thanks to this data neighbouring categories can be interpreted, for the same or for different variables. The axes show the main characteristics of the links existing between variables, and provide a synthetic description of our ten super-classes.

A supplementary variable was plotted, consisting of age divided into five subcategories : less than 25 years old (*m25*), from 25 to 35 (*25a35*), from 35 to 45 (*35a45*), from 45 to 55 (*45a55*) and over 55 years old (*p55*).

The graphs corresponding to the first three axes for the two regions studied can be found in the appendix II. Only our findings for the Nord region will be commented on here as those obtained for the second region differ very slightly.

The first axis relates essentially to the category of those doing occasional work. It reveals a sharp contrast between various measures of this work (*arp117*, *arm117*, *arm78*, *p3par*, *13par*) and the absence of occasional work over the whole period

(*0par*, *0ar*, *0peran*). The second axis is defined in terms of education, level of skill and the length of the latest period of unemployment. On the right of this graph there is clear evidence of a correlation between high skills and a high educational level on the one hand, and short duration of unemployment and the fact of an individuals first search for work on the other. On the left of the graph, very low skills (*mano*) and a poor educational level (*form0*) are correlated with a lengthy duration of unemployment (over 12 months). If we consider the supplementary variable, the second axis shows a major difference between the situation of young people (on the right) and older ones (on the left).

It should be noted that involvement in occasional work during unemployment is clearly orthogonal to skills, education, age and length of unemployment. On the basis of the data available jobless people in this category seem very close to the others (except for the fact that the latter do occasional work). These findings should however be modulated concerning the duration of unemployment by taking into account the fact that different groups may be formed within the category of those involved in occasional work, as will be seen from our ensuing analysis.

Though the third axis is not as well defined as the preceding ones, it does show the effect of a recurrence of unemployment which is not necessarily linked to the duration of the latest period. The fourth axis is defined by age whereas the fifth is defined by skill and not by educational level.

The ten super-classes obtained with the Kohonen classification may be defined as follows:

- *Group 1* consists of unemployed workers occupying occasional jobs, and who have a high level of education and skills. These main characteristics are also to be found in group 6.

- *Group 2* mainly comprises poorly educated individuals, aged over 35, where average duration of unemployment exceeded 12 months during their latest jobless period. Some of them are involved in occasional work.

- *Groups 3 and 4* are very close and are characterised by older poorly educated individuals, who do no occasional work. *Group 3* is slightly different in so far as a significant number of individuals here are trying to return to the labour market.

- *Groups 5, 9 and* 10 are composed of individuals who do occasional work, but who differ in terms of educational and skill level which is higher in *group 9* than in *group 5* and higher in *group 10* than in *group 9*, though the duration of their unemployment is constantly decreasing. One major difference characterising *group 10* is that these individuals appear to be often involved in occasional work (at least twice during each period of unemployment). It is also worth noting that the latter category experiences longer periods of unemployment.

- *Group 7* is rather different from the rest of the population studied as these workers undergo at least three periods of unemployment and most of them do occasional work. This clearly demonstrates a mechanism we explained earlier, involving the link between, fixed-term contracts, on the one hand,  and unemployment and occasional work on the other, the latter being presumably very similar to that done when working full-time;

- *Group 8* is characterised by a very long duration of unemployment. These individuals are older and do no occasional work.

These findings could be defined more precisely on the basis of further study of the different stages these groups went through over the three year period of observation. Overall, it can be observed that among those individuals working in occasional jobs there is a clear distinction between those who change jobs quite frequently and those for whom the occasional job becomes a permanent one.

**5. Factors influencing the classification**

We have been able to determine the most important variables used to construct the classification and the particular mix explaining each group obtained. This was done by applying a canonical discriminant analysis to the variable identifying the ten groups and the eleven quantitative variables used in the Kohonen classification.

A number of significant canonical components for the Nord region (as well as for the Rhone region) were obtained. The first four (representing 79 % and 81 % of the sum of the eigen values) are easy to interpret and support our earlier findings.

The first one consists of different variables relating to occasional work, and exerts a positive influence. This is the principal factor leading to the establishment of distinct classes. The second one is determined (for the Nord region) on the basis of age, length and recurrence of unemployment. The third one represents a brief period of unemployment without recurrence. The fourth one refers to the recurrence of unemployment and the amount of occasional work done over a given period of unemployment. The fifth canonical component reflects the number of job offers received per month of unemployment.

On the basis of the means computed for each class using our canonical variables we are now able to define a broad classification of the main categories of the unemployed. The first category consists of young people experiencing a short period of joblessness and who do no occasional work (*groups 1, 2, 6 and 7*). The second is characterised by older people undergoing a lengthy period of unemployment and who do no occasional work (*groups 3, 4 and 8*). The third comprises individuals who often resort to occasional work, who stay unemployed and experience no recurrence of unemployment (*group 5*). The fourth one represents those workers who do a lot of occasional work, experience both relatively short periods of unemployment and regular recurrence of it (*group 9 an 10*). The last two groups would seem to be involved in a permanent precarious situation. They are either employed for a short period and then become jobless (and perhaps resort to occasional work in the same field) or are generally registered as unemployed but work significantly fewer monthly hours than in a regular job.

These results are not contradictory with those obtained by Joutard (1992) on similar data, precisely on the question of the link between the total job seniority before the spell of unemployment and the greater opportunity to find a stable job. In fact, our results suggest that the practice of occasional jobs and the total duration of this practice is not perceived by future employers as a kind of accumulated human capital but as a negative signal of worker's potential. The consequence is that it exists a mechanism, very close in some way to the one enlightened by Joutard, between very occasional jobs during a spell of unemployment and a strong opportunity to find a precarious job. The recurrence of unemployment is a constitutive part of this mechanism and it illustrate the fact that the secondary segment is not a collection of heterogeneous activities produced by exits from the primary segment, but a coherent set of situations, possessing an internal logic and whose relations to the primary segment have to be explained. One of the results obtained here is that a longer duration of unemployment without practicing an occasional job seems to be the difficult way to return to the primary segment.

With these precisions, it is noticeable that comparable results are obtained by Joutard with 1986-87 data and over the more recent period we studied, 1993-96.

## 6. Conclusion

Quantitative classification enables classes of unemployed workers to be clearly differentiated. These classes may be easily defined by the methods we used. Some of the variables used are essential to building a typology of jobless workers. This typology is built around the variables relating to occasional work coupled with varying periods of unemployment and some degree of recurrence as well as around the age factor variables.

Granier (1998), using the whole population of the unemployed as if it was homogeneous, obtains biased results, mixing the different types of occasional jobs and, mainly, ignoring the segmentation of the labor market.

Our study provides only part of the evidence concerning the existence of segments in the labour market depending on entry-exit conditions. Much further investigation is required into the situation of the different classes of individuals who find themselves in and out of employment.

# APPENDIX I

The Kohonen algorithm (Kohonen, 1984, 93, 95; Cottrell, Fort, Pagès, 1997) is a well-known unsupervised learning algorithm which produces a map composed by a fixed number of units arranged across a network, generally a two-dimensional grid, with $n$ by $n$ units, but the method can be used with any topological organization of the Kohonen network.

A physical neighborhood relation between the units is defined and for each unit $i$, $V_r(i)$ represents the neighborhood with radius $r$ centered at $i$.

Each unit is characterized by a code vector $C_i$ of the same dimension as the input space (the data space).

The learning algorithm takes the following form :
  * at time 0, $C_i(0)$ is randomly defined for each unit $i$;
  * at time $t$, we present a data vector $x(t)$ randomly chosen among the rows of the data matrix and we determine the winning unit $i^*$, which minimizes the Euclidean distance between $x(t)$ and $C_i(t)$;
  * we then modify the $C_i$ in order to move the code vectors of the winning unit $i^*$ and its physical neighbors towards $x(t)$, using the following relations :

$$C_i(t+1) = C_i(t) + \varepsilon(t)\ (x(t) - C_i(t)) \text{ for } i \text{ in } V_{r(t)}(i^*) \quad (1)$$
$$C_i(t+1) = C_i(t) \text{ for other } i \quad (2)$$

where $\varepsilon(t)$ is a small positive adaptation parameter, $r(t)$ is the radius of $V_{r(t)}$ and $\varepsilon(t)$ and $r(t)$ are progressively decreased during the learning.

In fact the Kohonen algorithm in a generalization of the Forgy algorithm, where the update of the code vectors involve also the neighbors of the winning unit and not only the winning unit.

After learning, each unit $i$ is represented in the data space by its code vector $C_i$.

Then each observation is classified by a nearest neighbor method : observation $k$ belongs to class $i$ if and only if the code vector $C_i$ is the closest to the data vector $k$ among all the code vectors. The distance is the Euclidean distance in general, but it can be chosen in another way according to the context.

With respect to any other classification method, the main characteristic of the Kohonen classification is the conservation of the topology: after learning, « close » observations are associated to the same class or to « close » classes according to the definition of the neighborhood in the Kohonen network. This feature allows to consider the resulting classification as a good starting point to further developments and valuable representations as stated below.

As mentioned above, the first and raw result we get after learning, is a classification of the $N$ observations into $P = n \times n$ classes. Eventually some classes can be empty. At the same time, we get the code vectors (associated to each unit). The first two figures in Appendix II shows standard representations of the code vectors inside their own unit.

The choice of the number $P$ of units is arbitrary, and there does not exist any method to better choose the size of the network. We can guess that the « relevant » number of classes could often be smaller than $P$ (which is commonly equal to 100). It is also difficult to give relevant interpretation of a too large number of classes. So we propose to reduce the number of classes by means of a hierarchical classification of the $P$ code vectors using the Ward distance for example. As the code vectors are already organized across the grid, a standard method appears to be relevant.

In this way, we define two embedded classifications, and can distinguish the classes (Kohonen classes or « micro-classes ») and the « macro-classes » which group together some of the « micro-classes ». To make visible this two-levels classification, we affect to each « super-class » some colour or trame or grey level (here). Appendix II shows how the « micro-classes » are grouped together to constitute 10 « super-classes ».

The advantage of this double classification is the possibility to analyze the data set at a « super » level where general features emerge and at a « micro » level to determine the characteristics of more precise phenomena and especially the paths to go from one class to another one.

In the applications that we treated, the « super classes » create connected areas in the grid (one « super class » which contain some « micro class » contain also a neighbor of this class). This remark is very attractive because it confirms the topological properties of the Kohonen maps. Nevertheless, in some cases, it is possible to find a « super class » split into two pieces. In that case, one can guess that the data set is folded over and try to control it by studying the distortion.

# APPENDIX II

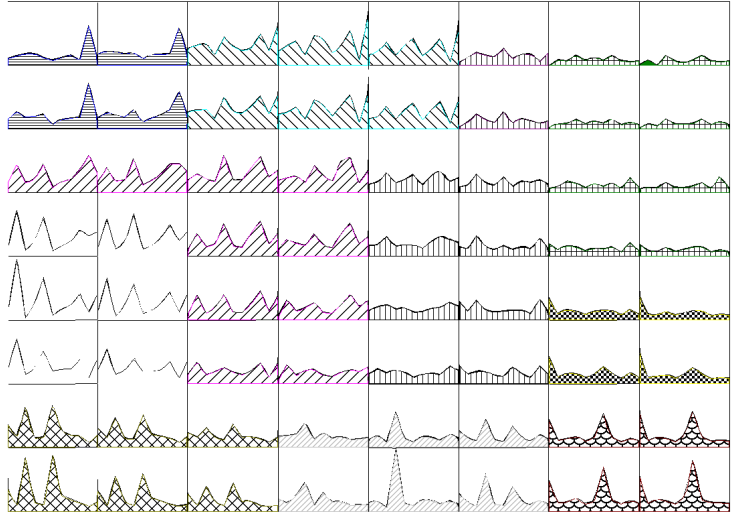

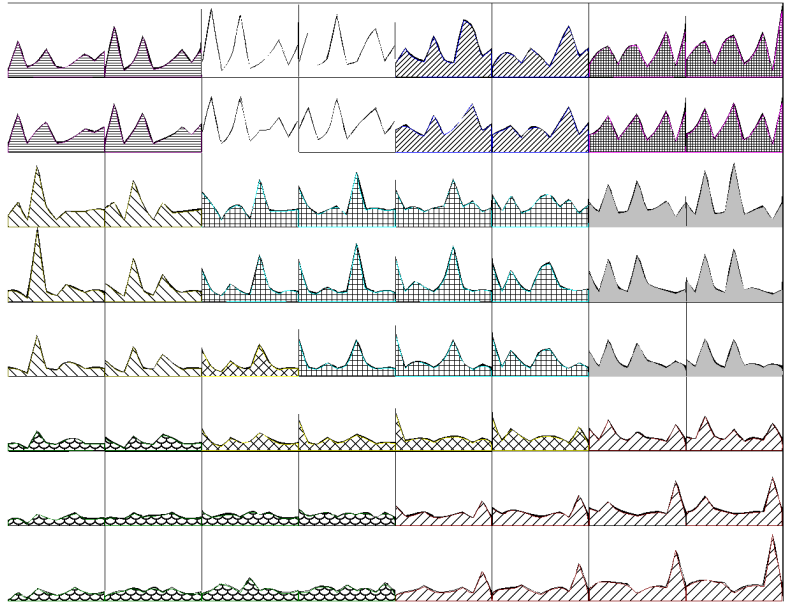

COMPOSITION OF THE GROUPS OBTAINED FROM KOHONEN CLASSIFICATION

| Gr. | Duration | | | Nb of unem. period | | | Category of age | | |
|---|---|---|---|---|---|---|---|---|---|
| | < 6 m. | 6-12 | > 12m | 1 d. | 2 d. | 3 & more | < 25 | 25-45 | > 45 |
| **NORD** | | | | | | | | | |
| 1 | 57.37 | 26.36 | 16.27 | 74.16 | 25.84 | . | 59.24 | 40.76 | 0.00 |
| 2 | 21.82 | 19.94 | 58.25 | 30.73 | 61.77 | 7.50 | 30.38 | 63.47 | 6.16 |
| 3 | 33.54 | 19.66 | 46.80 | 76.44 | 19.59 | 3.97 | 0.10 | 37.09 | 62.82 |
| 4 | 37.70 | 23.39 | 38.91 | 83.13 | 16.39 | 0.47 | . | 52.90 | 47.10 |
| 5 | 1.22 | 3.99 | 94.79 | 85.04 | 14.18 | 0.78 | 21.12 | 65.40 | 13.47 |
| 6 | 69.93 | 15.48 | 14.59 | 67.63 | 24.95 | 7.42 | 41.14 | 53.17 | 5.69 |
| 7 | 64.13 | 24.88 | 10.99 | . | . | 100.00 | 38.32 | 57.28 | 4.40 |
| 8 | 0.11 | 0.30 | 99.59 | 96.95 | 2.91 | 0.14 | 6.34 | 71.56 | 22.10 |
| 9 | 24.58 | 33.88 | 41.54 | 46.94 | 35.79 | 17.27 | 48.76 | 46.46 | 4.77 |
| 10 | 66.43 | 24.49 | 9.08 | 72.13 | 24.01 | 3.87 | 50.77 | 42.03 | 7.20 |
| **RHONE** | | | | | | | | | |
| 1 | 46.89 | 22.04 | 31.07 | 98.39 | 1.61 | . | 48.13 | 51.87 | 0.00 |
| 2 | 47.33 | 36.68 | 15.99 | 55.09 | 39.17 | 5.73 | 51.82 | 47.04 | 1.14 |
| 3 | 59.99 | 24.93 | 15.08 | . | 56.92 | 43.08 | 31.09 | 65.86 | 3.05 |
| 4 | 49.03 | 28.36 | 22.61 | 81.79 | 16.16 | 2.05 | 0.00 | 62.72 | 37.28 |
| 5 | 30.96 | 19.70 | 49.35 | 77.05 | 19.10 | 3.85 | 0.00 | 29.26 | 70.74 |
| 6 | 0.06 | 0.06 | 99.89 | 94.67 | 5.31 | 0.03 | 10.56 | 66.63 | 22.82 |
| 7 | 10.74 | 29.88 | 59.38 | 38.28 | 44.56 | 17.17 | 26.20 | 62.99 | 10.81 |
| 8 | 71.50 | 19.50 | 9.00 | 76.68 | 18.64 | 4.68 | 35.18 | 60.41 | 4.41 |
| 9 | 0.49 | 0.84 | 98.68 | 90.77 | 9.20 | 0.03 | 26.66 | 59.90 | 13.45 |
| 10 | 40.66 | 34.73 | 24.61 | 75.57 | 19.97 | 4.46 | 36.87 | 51.14 | 11.99 |

| Gr | Skills | | | | Educational level | | | | | |
|---|---|---|---|---|---|---|---|---|---|---|
| | sup | clerks | craft | na | sec +4 | sec. +2 | sec. dip | sec. | prim. | <prim. |
| **NORD** | | | | | | | | | | |
| 1 | 13.71 | 53.49 | 32.11 | 0.69 | 8.30 | 11.78 | 20.15 | 39.22 | 7.72 | 12.83 |
| 2 | 5.72 | 54.05 | 39.65 | 0.58 | 2.94 | 5.14 | 12.05 | 45.51 | 7.83 | 26.53 |
| 3 | 22.19 | 34.43 | 43.32 | 0.07 | 3.74 | 4.04 | 9.35 | 34.69 | 6.92 | 41.25 |
| 4 | 12.31 | 50.24 | 36.90 | 0.55 | 4.37 | 3.21 | 7.21 | 25.10 | 8.56 | 51.55 |
| 5 | 4.06 | 47.83 | 47.87 | 0.24 | 1.59 | 2.98 | 9.48 | 44.58 | 7.08 | 34.29 |
| 6 | 14.27 | 52.02 | 33.08 | 0.64 | 6.21 | 17.08 | 17.79 | 43.63 | 5.57 | 9.72 |
| 7 | 5.38 | 44.41 | 49.98 | 0.22 | 2.61 | 5.89 | 13.63 | 51.16 | 6.88 | 19.82 |
| 8 | 4.30 | 47.53 | 47.66 | 0.52 | 1.74 | 2.56 | 5.44 | 30.61 | 7.37 | 52.28 |
| 9 | 10.91 | 49.50 | 39.09 | 0.50 | 4.86 | 10.35 | 20.04 | 45.57 | 5.16 | 14.01 |
| 10 | 17.44 | 49.32 | 32.75 | 0.48 | 9.23 | 17.34 | 21.40 | 35.51 | 4.49 | 12.03 |

| | | | | | RHONE | | | | | |
|---|---|---|---|---|---|---|---|---|---|---|
| 1 | 17.94 | 48.64 | 23.35 | 10.08 | 13.45 | 15.40 | 16.95 | 37.27 | 5.75 | 11.19 |
| 2 | 17.93 | 16.10 | 45.76 | 24.29 | 13.85 | 10.51 | 16.06 | 19.59 | 41.4 | 4.74 |
| 3 | 10.06 | 46.34 | 29.35 | 14.25 | 6.19 | 10.70 | 14.66 | 46.02 | 7.49 | 14.94 |
| 4 | 16.89 | 47.25 | 25.75 | 10.10 | 9.22 | 7.95 | 12.21 | 32.44 | 6.69 | 31.50 |
| 5 | 25.97 | 38.43 | 35.18 | 0.42 | 7.74 | 6.39 | 9.36 | 34.21 | 6.08 | 36.22 |
| 6 | 8.06 | 42.62 | 29.40 | 19.92 | 5.76 | 5.93 | 10.33 | 34.59 | 7.29 | 36.09 |
| 7 | 12.00 | 40.95 | 28.13 | 18.91 | 7.92 | 10.63 | 14.16 | 45.26 | 5.95 | 16.09 |
| 8 | 18.68 | 46.18 | 34.36 | 0.77 | 6.95 | 16.32 | 14.91 | 47.23 | 5.09 | 9.50 |
| 9 | 7.91 | 39.72 | 23.55 | 28.82 | 5.75 | 7.25 | 14.08 | 43.90 | 6.66 | 22.37 |
| 10 | 18.95 | 46.90 | 27.86 | 6.29 | 9.27 | 15.65 | 19.30 | 39.14 | 4.37 | 12.26 |

| Gr | Reasons of job loss | | | | | Type of unemployment exit | | | | |
|---|---|---|---|---|---|---|---|---|---|---|
| | end | lay-off | 1st job | return | na | job | p.job | quit | train. | susp. |
| | | | | | NORD | | | | | |
| 1 | 42.91 | 9.36 | 37.44 | 8.73 | 1.56 | 24.54 | 19.40 | 13.90 | 5.33 | 9.92 |
| 2 | 62.25 | 8.73 | 12.82 | 9.90 | 6.30 | 17.49 | 11.25 | 10.54 | 9.98 | 18.41 |
| 3 | 48.21 | 41.81 | . | 8.40 | 1.58 | 21.17 | 7.48 | 10.70 | 5.55 | 21.33 |
| 4 | 48.63 | 25.21 | 5.46 | 18.87 | 1.82 | 16.32 | 16.21 | 10.89 | 4.70 | 16.79 |
| 5 | 59.04 | 14.93 | 6.87 | 8.84 | 10.33 | 31.55 | 5.25 | 9.55 | 11.6 | 17.94 |
| 6 | 56.88 | 14.91 | 17.47 | 9.66 | 1.09 | 41.97 | 17.59 | 5.89 | 4.61 | 5.50 |
| 7 | 83.37 | 3.98 | 5.42 | 4.62 | 2.61 | 31.38 | 14.08 | 7.87 | 5.45 | 9.30 |
| 8 | 42.66 | 25.56 | 7.26 | 23.08 | 1.44 | 8.59 | 7.78 | 9.43 | 11.0 | 18.65 |
| 9 | 63.82 | 9.76 | 16.20 | 6.23 | 3.99 | 38.52 | 12.04 | 8.94 | 6.18 | 10.04 |
| 10 | 55.12 | 13.09 | 23.04 | 6.76 | 1.98 | 43.72 | 13.00 | 8.31 | 1.26 | 5.46 |
| | | | | | RHONE | | | | | |
| 1 | 40.25 | 17.29 | 24.40 | 15.15 | 2.92 | 29.82 | 26.00 | 13.02 | 7.41 | 9.27 |
| 2 | 52.03 | 13.68 | 18.02 | 11.00 | 5.27 | 35.40 | 19.52 | 9.44 | 4.49 | 8.87 |
| 3 | 60.61 | 12.39 | 6.57 | 9.99 | 10.4 | 25.77 | 19.15 | 8.07 | 6.10 | 17.39 |
| 4 | 38.06 | 35.71 | 3.77 | 17.68 | 4.79 | 26.03 | 19.63 | 11.40 | 5.01 | 15.27 |
| 5 | 35.46 | 51.35 | 0.04 | 9.34 | 3.81 | 22.45 | 13.15 | 9.51 | 1.60 | 24.64 |
| 6 | 38.99 | 32.43 | 6.87 | 17.17 | 4.54 | 18.50 | 11.78 | 9.39 | 15.0 | 18.36 |
| 7 | 55.89 | 19.18 | 6.69 | 8.44 | 9.81 | 31.44 | 14.49 | 6.24 | 7.69 | 15.72 |
| 8 | 49.64 | 22.23 | 12.00 | 15.45 | 0.68 | 33.64 | 25.73 | 9.09 | 1.55 | 5.05 |
| 9 | 48.33 | 23.07 | 8.92 | 11.01 | 8.68 | 35.37 | 9.23 | 6.31 | 12.8 | 12.26 |
| 10 | 54.39 | 20.15 | 13.51 | 10.08 | 1.87 | 27.82 | 16.09 | 6.60 | 1.78 | 5.62 |

| Gr | Monthly hours in occ. jobs | | | | | Nb of occ.job by year of unemp. | | | |
|---|---|---|---|---|---|---|---|---|---|
| | 0 | 0-39 | 39-78 | 78-117 | >117 | 0 | 0-0.5 | 0.5-2 | >2 |
| **NORD** | | | | | | | | | |
| 1 | 99.97 | 0.03 | . | . | . | 99.39 | 0.31 | 0.31 | . |
| 2 | 77.08 | 3.80 | 9.52 | 3.51 | 6.10 | 52.00 | 24.70 | 23.24 | 0.06 |
| 3 | 86.02 | 1.74 | 5.58 | 1.35 | 5.32 | 81.16 | 6.27 | 11.75 | 0.82 |
| 4 | 97.74 | 0.66 | 1.24 | 0.18 | 0.18 | 96.68 | 0.87 | 2.44 | . |
| 5 | 0.24 | 1.02 | 7.38 | 11.41 | 79.96 | . | 25.96 | 67.50 | 6.53 |
| 6 | 86.50 | 2.56 | 5.25 | 2.11 | 3.58 | 79.27 | 1.34 | 15.42 | 3.97 |
| 7 | 78.34 | 3.73 | 6.37 | 2.55 | 9.02 | 40.55 | 13.76 | 43.68 | 2.01 |
| 8 | 83.99 | 2.80 | 7.56 | 1.77 | 3.89 | 83.77 | 15.58 | 0.65 | . |
| 9 | 6.38 | 5.25 | 15.44 | 10.13 | 62.80 | . | 8.33 | 81.39 | 10.28 |
| 10 | 3.72 | 9.95 | 25.80 | 12.66 | 47.87 | . | . | 20.10 | 79.90 |
| **RHONE** | | | | | | | | | |
| 1 | 97.31 | 1.13 | 1.45 | 0.11 | . | 97.31 | 0.75 | 1.93 | . |
| 2 | 11.43 | 12.39 | 47.79 | 19.05 | 9.33 | . | 0.21 | 68.63 | 31.16 |
| 3 | 80.95 | 3.07 | 11.15 | 3.10 | 1.73 | 62.29 | 12.49 | 24.68 | 0.53 |
| 4 | 96.68 | 0.84 | 2.15 | 0.33 | . | 94.28 | 1.90 | 3.82 | . |
| 5 | 80.88 | 2.50 | 8.16 | 4.87 | 3.58 | 77.69 | 6.22 | 14.73 | 1.37 |
| 6 | 76.53 | 3.69 | 9.59 | 8.60 | 1.59 | 76.45 | 21.11 | 2.44 | . |
| 7 | . | 0.22 | 1.75 | 16.35 | 81.68 | . | 19.32 | 79.00 | 1.67 |
| 8 | 88.09 | 2.18 | 6.05 | 1.86 | 1.82 | 84.95 | 0.45 | 11.82 | 2.77 |
| 9 | 0.03 | 0.87 | 10.63 | 30.66 | 57.80 | . | 8.12 | 84.08 | 7.80 |
| 10 | 0.36 | 2.50 | 20.24 | 16.90 | 60.01 | . | 0.40 | 44.09 | 55.51 |

| Gr | Nb of occ.job by period of unemp. | | | Prop. of occ. jobs in unem. duration | | | |
|---|---|---|---|---|---|---|---|
|  | 0 | 0-1 | >1 | 0 | 0-0.1 | 0.1-0.3 | >0.3 |
| **NORD** | | | | | | | |
| 1 | 99.39 | 0.61 | . | 99.39 | 0.61 | . | . |
| 2 | 52.00 | 45.99 | 2.00 | 52.00 | 32.78 | 12.28 | 2.94 |
| 3 | 81.16 | 17.23 | 1.61 | 81.16 | 10.27 | 6.89 | 1.67 |
| 4 | 96.68 | 3.32 | . | 96.68 | 2.44 | 0.80 | 0.07 |
| 5 | . | 21.43 | 78.57 | . | 13.34 | 36.22 | 50.44 |
| 6 | 79.27 | 19.39 | 1.34 | 79.27 | 7.49 | 10.24 | 3.01 |
| 7 | 40.55 | 57.76 | 1.69 | 40.55 | 28.93 | 26.76 | 3.76 |
| 8 | 83.77 | 14.36 | 1.88 | 83.77 | 12.56 | 2.75 | 0.92 |
| 9 | . | 79.31 | 20.69 | . | 14.19 | 44.23 | 41.58 |
| 10 | . | 66.18 | 33.82 | . | . | 17.58 | 82.42 |
| **RHONE** | | | | | | | |
| 1 | 97.31 | 2.69 | . | 97.31 | 2.44 | 0.24 | . |
| 2 | . | 79.02 | 20.98 | . | 12.07 | 61.65 | 26.28 |
| 3 | 62.29 | 36.04 | 1.67 | 62.29 | 23.38 | 12.67 | 1.66 |
| 4 | 94.28 | 5.60 | 0.13 | 94.28 | 4.10 | 1.47 | 0.15 |
| 5 | 77.69 | 18.08 | 4.24 | 77.69 | 10.45 | 9.80 | 2.06 |
| 6 | 76.45 | 18.81 | 4.74 | 76.45 | 18.98 | 3.92 | 0.65 |
| 7 | . | 76.96 | 23.04 | . | 24.04 | 49.91 | 26.05 |
| 8 | 84.95 | 13.45 | 1.59 | 84.95 | 4.68 | 8.27 | 2.09 |
| 9 | . | 5.44 | 94.56 | . | 12.96 | 50.73 | 36.31 |
| 10 | . | 58.40 | 41.60 | . | . | 3.03 | 96.97 |

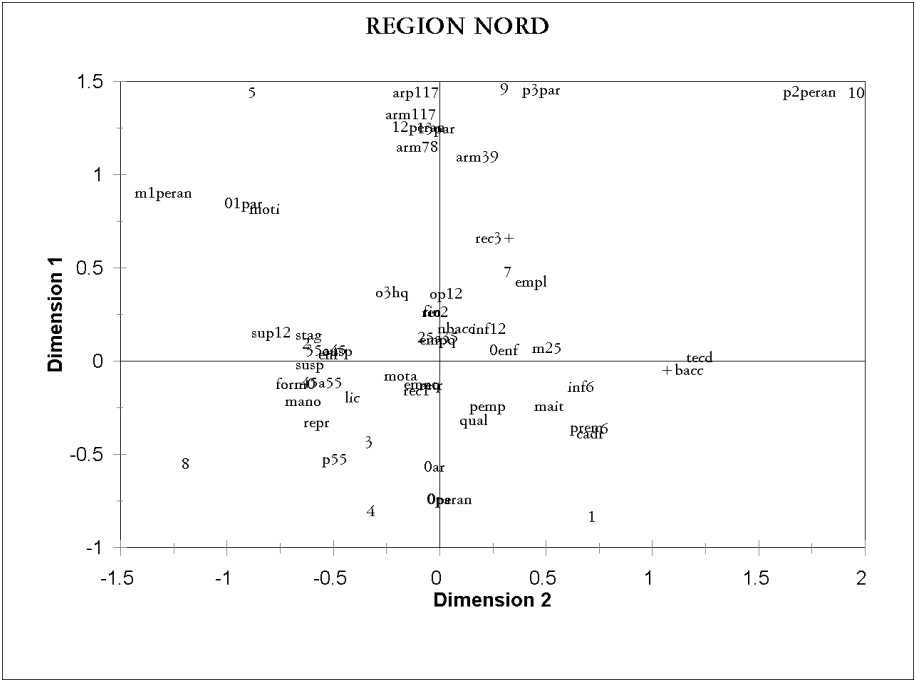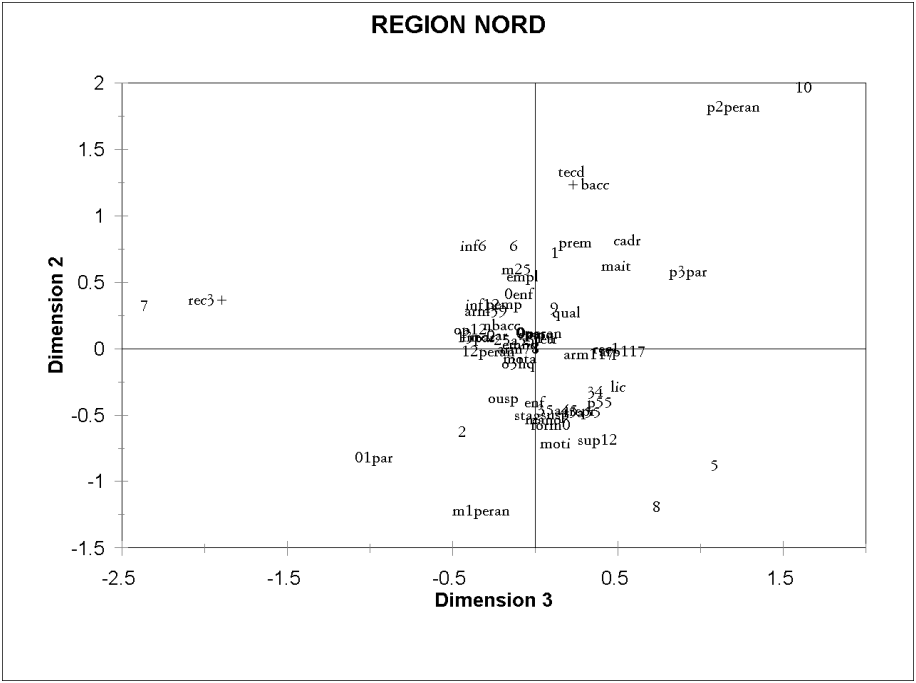

## REGION RHONE

(Scatter plot: Dimension 1 vs Dimension 2)

## REGION RHONE

(Scatter plot: Dimension 2 vs Dimension 3)